\documentstyle[12pt]{article}

\begin{document}
\begin{titlepage}
\rightline{Preprint  ISI/1997/11}
\rightline{November 1997}
\rightline{q-alg/9711xxx}
\vspace{1.5cm}

\begin{center} 
{\huge \bf Twist Deformation of the }\\[.5cm]
{\huge \bf rank one Lie Superalgebra}
\end{center}
\vspace{.6cm}

\centerline{\Large \bf E.Celeghini}
\vspace{.3cm}
\begin{center}
{Dipartimento di Fisica del'Universit\`a and INFN, Firenze}\\
{L.go Fermi 2 - I50125 Firenze, Italy}
\end{center}

\vspace{.3cm}
\centerline{\bf and}
\vspace{.3cm}

\centerline{\Large \bf  P.P.Kulish
\footnote{On leave
of absence from the St.Petersburg 
Department of the Steklov
Mathematical Institute, 
Fontanka 27, St.Petersburg, 191011,
Russia \,; ( kulish@pdmi.ras.ru )} }
\vspace{.3cm}
\begin{center}
{Institute for Scientific Interchange Foundation}\\
{Villa Gualino, 10133, Torino, Italy}
\end{center}

\vspace{1.0cm}

\begin{abstract}
The Drinfeld twist is applyed to deforme the rank one 
orthosymplectic Lie superalgebra $osp(1|2)$. The twist 
element is the same as for the $sl(2)$ Lie algebra due to 
the embedding of the $sl(2)$ into the superalgebra $osp(1|2)$. 
The $R$-matrix 
has the direct sum structure in the irreducible representations 
of $osp(1|2)$. The dual quantum group is defined using the 
FRT-formalism. It includes the Jordanian quantum group $SL_\xi(2)$ 
as subalgebra and Grassmann generators as well. 
\end{abstract}

\end{titlepage}

\section{The deformed algebra $osp_{\xi}(1|2)$}

It is difficult to overestimate the role of the rank one 
Lie algebra $sl(2)$
in the theory of Lie groups and their applications.
The corresponding role for Lie superalgebras is played by the 
orthosymplectic superalgebra $osp(1|2)$ with
five generators $\{h, X_-, X_+, v_-, v_+\}$ and commutation relations
(Lie super- or ${\bf Z}_2$ graded-brackets): 

\begin{eqnarray}
[ h, X_{\pm} ] = \pm 2 X_{\pm}, \;\;\;\; 
[ X_{+}, X_{-} ] = h \;,  \label{one} 
\end{eqnarray}
\begin{eqnarray}
[ h, v_{\pm} ] = \pm v_{\pm} ,  \;\; \;\;
[ v_+, v_- ]_+ = - h / 4 \;, \label{two}   
\end{eqnarray}
\begin{eqnarray}
[ X_{\pm}, v_{\pm} ] = 0, \;\; \;\;
[ X_{\pm} , v_{\mp} ] = v_{\pm}, \;\;\;\; 
[ v_{\pm}, v_{\pm}]_+ = \pm X_{\pm} / 2  \;. \label{three} 
\end{eqnarray}

\noindent The generators $h$ and $ X_{\pm}$ are even (zero parity $p=0$), 
while $v_{\pm} $ are odd, $p = 1$. 
As a Hopf superalgebra, the universal enveloping ${\cal U}\bigl(osp(1|2)\bigr)$ 
of $osp(1|2)$ is 
generated, as $sl(2)$, just by three elements: it is sufficient to start from
$\{h, v_-, v_+\}$ restricted by the relations (\ref{two}) only, and define 
$X_{\pm} \equiv \pm 4 v_{\pm}^2$. 

The quantum deformation of $sl(2)$ can be considered as a "pivot" 
of the quantum group theory \cite{Dr,FRT}, while
the corresponding quantum superalgebra $osp_q(1|2)$ constructed 
in \cite{K,KR,S}, is the corresponding analogue 
for the quantum supergroups.
As a quasitriangular Hopf superalgebra $osp_q(1|2)$, analogously to 
the universal enveloping of $osp(1|2)$,
is generated by  three elements $\{h, v_-, v_+\}$ under the relations 
\begin{eqnarray}
[ h, v_{\pm} ] = \, \pm \; v_{\pm} ,\hskip 1.5cm [ v_+, v_-] = \, - \;
\frac{1}{4} \,(q^h - q^{-h})/ (q - q^{-1})  .\nonumber
\end{eqnarray}
\noindent 
It is worthy to note that, while $sl(2)$ is embedded into $osp(1|2)$,
such embedding does not exist for $sl_q(2)$ into $osp_q(1|2)$
because the coproduct of even elements $X_{\pm} \sim v_{\pm}^2$ 
includes also odd ones.

The aim of this paper is to construct and study 
the twist deformation \cite{Drtw} of $osp(1|2)$ 
that looks, in some sense,
more natural than $osp_q(1|2)$ because it is consistent 
with this fundamental property of 
inclusion $ sl(2) \subset osp(1|2)$ 
and it is generated by the same twist element of $sl(2)$.

The triangular Hopf algebra $sl_\xi(2)$ 
(cf. \cite{Z,Og,Vl,KhST,BH,ADM}, and Refs therein) 
is given by the extension of the twist deformation of 
the universal enveloping of the Borel sub-algebra $B_- \equiv
\{ h, X_-\}$ to the whole ${\cal U}\bigl(sl(2)\bigr)$. The twist element 
${\cal F}$ is
\begin{eqnarray} 
{\cal F} =\; 1\; +\;\; \xi \; h \otimes X_-\;\; +\;\; \frac{\xi^2}{2} \; h(h+2) 
\otimes X_-^2 \; +\; \dots \nonumber 
\end{eqnarray}
that can be written as
\begin{eqnarray}
{\cal F}\; =\; \bigl( 1 - 2 \,\xi \; 1 \otimes 
X_-\bigr)^{-\frac{1}{2}(h\otimes 1)}
\; =\; {\rm exp} \bigl(\frac{1}{2} \; h \otimes \sigma \bigr) 
\label{f}
\end{eqnarray} 
where $\sigma = \, - {\sl ln} \, (1-2 \, \xi \, X_-)$.

Let us recall from \cite{Drtw} that for a quasitriangular Hopf 
algebra ${\cal A}$ with an $R$-matrix ${\cal R}$ the twisted Hopf 
algebra ${\cal A}_t$ has $R$-matrix ${\cal R}^{({\cal F})}$ 
given by the twist transformation  
\begin{equation} 
{\cal R}^{({\cal F})}={\cal F}_{21}{\cal RF}^{-1}  
\label{a2.2} 
\end{equation}
of the original $R$-matrix ${\cal R},$ where ${\cal 
F}_{21}={\cal PFP},$ and ${\cal P}$ is the permutation map in 
${\cal A}\otimes {\cal A}$. The algebraic sector of ${\cal A}_t$ is 
not changed, and new coproduct is 
$\Delta_t = {\cal F} \Delta {\cal F}^{-1}$. 
The twist element satisfies the 
relations  in ${\cal A}\otimes {\cal A}$ \cite{Drtw} 
$$ 
(\epsilon \otimes id)\, {\cal F} = (id \otimes \epsilon)\, {\cal F} = 1\,, 
$$ 
and in ${\cal A}\otimes {\cal A}\otimes {\cal A}$ 
$$
{\cal F}_{12}\,(\Delta \otimes id)\, {\cal F} = 
{\cal F}_{23}\, (id \otimes \Delta) \,{\cal F} \,. 
$$

According to this Drinfeld definition, 
the algebraic relations of eqs. (\ref{one}) 
for the twisted $sl(2)$ are still the same, while the twisted coproduct 
$\Delta_t \equiv {\cal F} \Delta {\cal F}^{-1}$ is now on the 
generators 
\begin{eqnarray}
\Delta_t (h) &=& h \otimes e^\sigma + 1 \otimes h \;, \nonumber \\
\Delta_t (X_-) &=& X_- \otimes 1 + 1 \otimes X_- - 2 \, \xi \; X_- \otimes X_- =
\; X_- \otimes e^{-\sigma} + 1 \otimes X_- \,, \nonumber \\
\Delta_t (X_+) &=& X_+ \otimes e^\sigma + 1 \otimes X_+ 
               - \xi \, h \otimes e^\sigma h 
            + \frac{\xi}{2} \, h(h-2) \otimes e^\sigma (1-e^\sigma) .\nonumber 
\end{eqnarray}
\noindent 
Let us stress that this twist of the whole $sl(2)$ is obtained 
due to the embedding $B_- \subset sl(2)$.

Thus, knowing that $B_- \subset sl(2) \subset osp(1|2)$, 
the procedure can be simply
iterated to find $osp_\xi(1|2)$ (as well as the twisted deformations of all 
others nontrivial embeddings of
$B_-$). It is an easy  exercise, keeping in mind the expression of ${\cal F}$ 
 $\bigl($eq. (\ref{f})$\bigr)$, commutation relations (2), (3) 
and the primitive coproduct of $osp(1|2)$, to obtain:  
\begin{eqnarray}
\Delta_t (h) &=& h \otimes e^\sigma + 1 \otimes h \;, \nonumber \\
\Delta_t (v_-) &=& v_- \otimes e^{-\sigma/2} + 
1 \otimes v_- \;, \label{def} \\
\Delta_t (v_+) &=& v_+ \otimes e^{\sigma/2} 
                 + 1 \otimes v_+ +  \xi h \otimes v_- e^{\sigma} .  \nonumber
\end{eqnarray} 
One can reproduce the coproducts of $X_{\pm} $ by squaring the coproducts 
of $v_{\pm} $, taking into account the $Z_2$-grading of tensor product: 
$$ 
(x \otimes y)\, (u \otimes w) = (-1)^{p(u) p(y)} (x \,u\otimes y\,w) \;,  
$$ 
and the commutation relations (2), (3). 

The maps of counit $\epsilon $ and antipode $S$, necessary 
for a Hopf superalgebra definition, are 
\begin{eqnarray} 
\epsilon (h) &=& \epsilon (v_{\pm}) = 0 \,, \quad \epsilon (1) = 1 \,, \label{co}\\ 
S(h) &=& -h e^{-\sigma} \,, \quad S(v_-) = -v_- e^{\sigma /2} \,, \quad  
S(v_+) =  -(v_+ - \xi h v_- )e^{-\sigma /2} \,. \nonumber
\end{eqnarray}

We can thus arrive to the following \\ 
$\underline{Definition.}$ The Hopf superalgebra generated by three
elements 
$\{h,v_-,v_+\}$ satisfying the relations (\ref{two}), (\ref{def}) 
and (\ref{co}) is said to be
the twist deformation of  ${\cal U}\bigl(osp(1|2)\bigr)$ or $osp_\xi(1|2)$. 

This is a triangular Hopf superalgebra 
$( {\cal R}_{21}{\cal R} = 1)$ with universal $R$-matrix 
\begin{eqnarray} 
{\cal R} =\; {\rm exp} \bigl(\frac{1}{2} \; \sigma \otimes h \bigr)\; 
{\rm exp} \bigl(-\frac{1}{2} \; h \otimes \sigma \bigr) \,.   \label{R}
\end{eqnarray} 

The irreducible finite dimensional representations of $osp_\xi(1|2)$ 
$$ 
\rho_s : osp_\xi(1|2) \longrightarrow End \,(W_s) 
$$ 
are the same as for $osp(1|2)$, due to the unchanged algebraic 
relations (2). They are parametrized by the half-integer spin 
$s = 0, \frac {1}{2}, 1, \,...$, have dimension $4s + 1$, 
and are decomposed into 
a direct sum of two irreps of the $sl(2)$ \cite{Ritt}: 
$W_s = V_s + V_{s - \frac{1}{2}}$. Hence, the $R$-matrix in the irreps 
of $osp_\xi(1|2)$ is a direct sum of four $R$-matrices of $sl_\xi(2)$. 
For the first non-trivial case $s = 1/2$ one gets 
\begin{equation} 
{\bf R} = (\rho_{\frac {1}{2}} \otimes \rho_{\frac {1}{2}}) \,{\cal R} = 
R(\xi) + I_2 + I_2 + 1 \,,    \label{R1/2}
\end{equation} 
where $I_2$ are $2 \times 2$ unit matrices, and 
$ R(\xi)$ is the Jordanian solution 
to the Yang-Baxter equation (cf. \cite{Z}) 
\begin{equation}
R(\xi)=\left(
\begin{array}{ccccccccc}
1 & 0 & 0 & 0 \\
- \xi & 1 & 0 & 0 \\
\xi & 0 & 1 & 0 \\
\xi^2 & - \xi & \xi & 1 \\
\end{array}
\right) \,.  \label{Rxi}
\end{equation} 
The twist parameter can be scaled: $ \xi \rightarrow \exp (2u) \, \xi$ 
by the similarity transformation with the element $ \exp (-u\,h)$. 

The basis of the irreps tensor product decomposition 
will include deformed 
Clebsch-Gordan coefficients, expressed as linear combinations 
of the usual ones and the matrix elements of the 
twist $\cal F$ \cite{KS}. This is reflected in the spectral decomposition 
of the $R$-matrix itself in the tensor product $W_s \otimes W_l$ 
$$ 
{\hat R}^{s,l} = F^{s,l} \, 
(\sum_{j=|s-l|}^{s+l} ({\pm}) P^{j} )  \,(F^{s,l})^{-1} \;,
$$ 
where  $P^{j}$ are projectors onto irreducible representations of $osp(1|2)$.

\section{Quantum supergroup $OSp_\xi(1|2)$} 

The self-dual character of the twisted Borel subalgebra $(B_-)_\xi$ 
was pointed out
in~\cite{Og}. This is obvious in terms of the generators 
$\{h,\sigma\} \in(B_-)_\xi$
and the generators $\{s, p \} \in (B_-)_\xi'$~ of the dual,  
with the only non-trivial evaluations 
$\langle h, s \rangle = 2 \, ,\, \langle \sigma, p \rangle = 2 $ 
\cite{Og,Vl}:
\begin{eqnarray} \nonumber 
[ h , \sigma ] &=& 2 \, \bigl( 1 - e^{\sigma} \bigr)\,, \qquad 
[ p , s ] = 2 \, \bigl( 1 - e^{s} \bigr) \,, \\ 
\Delta (\sigma) &=& \sigma \otimes 1 + 1 \otimes \sigma \; , \hskip 1cm
\Delta (s) = s \otimes 1 + 1 \otimes s \;, \nonumber \\ 
\Delta (h) &=& h \otimes e^{\sigma} + 1 \otimes h \;,  \hskip 1cm 
\Delta (p) = p \otimes e^{s} + 1 \otimes p \;, \nonumber \\  
\epsilon (h) &=& \epsilon(\sigma ) = 0 \;, \qquad 
\epsilon (s) = \epsilon (p) = 0 \;, \nonumber \\  
S(h) &=& - h e^{-\sigma}\;, \; 
S(\sigma) = - \sigma \; , \qquad  
S(p) = - p e^{-s}\;, \; S(s) = - s \;. \nonumber 
\end{eqnarray}

The situation is different for the twisted Hopf super-subalgebra $(sB_-)_\xi$.
The latter is generated by two elements $\{ h , v_- \}$ as $(B_-)_\xi$ . 
However, due to the $Z_2$-grading its basis as a linear space 
consists of even $ \sigma^m h^n$ and odd $ \sigma^m v_- h^n$ elements 
($ \sigma =  - {\sl ln} \, (1 + 8 \, \xi \, v_-^2)$). \\
$\underline{Proposition}.$ The dual $(sB_-)_\xi'$ 
of the twisted Hopf superalgebra $(sB_-)_\xi$ is 
generated by three elements $\{ \nu , \eta , x \}$ satisfying the relations 
\begin{eqnarray} 
[\nu , \eta] &=& 0 \;, \quad 
[\nu , x] = \frac{1}{2} \bigl( 1 - e^{-2\nu} \bigr)\, \;, 
\quad [x , \eta] = \frac{1}{2} \eta \;, \quad  \eta^2 = 0 \;, \\
\Delta(\nu) &=& \nu \otimes 1 + 1 \otimes \nu \; , \hskip 1cm
\Delta(\eta) = \eta \otimes 1 + e^{-\nu} \otimes \eta \;, \nonumber \\
\Delta(x) &=& x \otimes 1 + e^{-2\nu} \otimes x + 
\frac{1}{8\xi} e^{-\nu } \eta \otimes \eta \;, \nonumber \\ 
\epsilon (x) &=& \epsilon(\eta ) = \epsilon (\nu) = 0 \;, \nonumber \\  
S(\eta) &=& - \eta e^{\nu}\;, \quad  
S(\nu) = - \nu \; , \quad  
S(x) = - x e^{2\nu} \;. \nonumber 
\end{eqnarray} 

One can check this by a straightforward calculation of evaluating the 
dual basis $ x^k \eta^\delta \nu^l$ of $(sB_-)_\xi'$ and 
$ \sigma^m v_-^\delta h^n$ of $(sB_-)_\xi \,,\; k, l, m, n = 0, 1, 2, ...; 
\delta = 0, 1$ with the only non-zero evaluations among the 
generators: $\langle h, \nu \rangle = 1 \, ,\, \langle v_-, \eta \rangle = 1
\;, \, \langle \sigma, x \rangle = 1 $. 
We shall prove it below by a reduction from the quantum 
supergroup $OSp_\xi(1|2)$. The universal $T$-matrix 
(bicharacter) is given in term of these basis by a product of three 
exponents 
$$ 
{\cal T} = \exp ({\sigma \otimes x}) \, \exp ({v_- \otimes \eta}) \, 
\exp ({h \otimes \nu}) \;. 
$$ 
It is interesting to point out that starting from a Hopf superalgebra 
without nilpotent elements we were forced to introduce Grassmann 
variables ($\eta$) in the dual superalgebra. 

The dual of the twisted Hopf superalgebra $osp_\xi(1|2)$ can be introduced 
using a $Z_2$-graded version of 
the FRT-formalism~\cite{FRT}, because the $R$-matrix in the fundamental 
representation is known (\ref{R1/2}). The $T$-matrix of generators of 
quantum supergroup $OSp_\xi(1|2)$ in this representation 
has dimension $3 \times 3$. 
There are two convenient basis in 
this irrep as $C^{3}$: i) with grading $(0,\; 1,\; 0)$ and ii) 
with grading $(0,\; 0,\; 1)$. The odd generators $v_-,\; v_+$ of $osp(1|2)$ 
are lower and upper triangular in the former basis, while the latter one 
is more convenient to write $\bf T$ in a block matrix form. These forms are 
\begin{equation}
{\bf T}=\left(
\begin{array}{ccc}
a & \alpha & b \\
\gamma & g & \beta \\
c & \delta & d
\end{array}
\right) \;, \qquad  {\bf T}=\left(
\begin{array}{cc}
T & \psi \\
\omega & g
\end{array}
\right) \;, 
\label{T}
\end{equation} 
where $T$ is $2\times 2$ matrix of the even generators $\{ a, b, c, d \}$, while 
$\psi $ and $ \omega$ are two component column $(\alpha \,, \; \delta)^t$ 
and row $(\gamma \,, \; \beta)$ vectors of odd elements. 

The $3 \times 3$ matrix {\bf T} of the $OSp_\xi(1|2)$ generators 
satisfies the FRT-relation 
\begin{eqnarray}
{\bf R T_1 T_2 = T_2 T_1 R} 
\label{FRT} 
\end{eqnarray}
with   ${\bf Z}_2$-graded tensor product and $9 \times 9$ 
$R$-matrix {\bf R} (\ref{R1/2}). From the block-diagonal 
form of {\bf R} (\ref{R1/2}) it follows for $2 \times 2$ matrix $T$
\begin{eqnarray}
R(\xi) T_1 T_2 = T_2 T_1 R(\xi) \;. 
\label{SL} 
\end{eqnarray}
Hence, one reproduces the algebraic sector (commutation 
relations) of the twisted quantum group $SL_{\xi}$(2) for 
the generators  
$\{ a, b, c, d \}$~\cite{Z}. For the other blocks 
of different dimension we get from (\ref{FRT})
\begin{eqnarray}
R(\xi) T_1 \psi_2 &=& \psi_2 T_1 \;, \quad g \;{\bf T} =  {\bf T} \;g \;, \\
\omega_1 T_2 &=& T_2 \omega_1 R(\xi)\;, \quad 
\omega_1 \psi_2 = - \psi_2 \omega_1 \;, \\
\omega_1 \omega_2 &=& - \omega_2 \omega_1 R(\xi)\,, \; 
R(\xi) \psi_1 \psi_2 = - \psi_2 \psi_1 \;. 
\label{OSp} 
\end{eqnarray} 
From the relations (\ref{SL}) - (\ref{OSp}) one gets centrality 
of the following elements: 
$$ 
{\det }_\xi T = a(d - \xi b) -cb \;, \quad g \;, \quad 
\theta = \omega \, T^{-1} \,\psi \;. 
$$ 

Coproduct, counit and antipode are given by the 
standard expressions of the FRT-formalism \cite{FRT}
\begin{eqnarray}
\Delta({\bf T}) = {\bf T} \otimes {\bf T} \;, \quad 
\epsilon ({\bf T}) = I_3 \;, \qquad S({\bf T}) = {\bf T}^{-1} \;. 
\label{CO}
\end{eqnarray} 
The inverse of $ {\bf T}$ is expressed in terms of the generators 
(\ref{T}) provided invertability of 
$\det_\xi T  \,,$ and $( g  - \omega \, T^{-1} \,\psi ) $
\begin{equation}
{\bf T}^{-1}=\left(
\begin{array}{cc}
I_2 & - T^{-1} \psi \\
0 & 1
\end{array}
\right) \, 
\left(
\begin{array}{cc}
T^{-1} & 0 \\
0 & ( g - \theta )^{-1}
\end{array}
\right) \, 
\left(
\begin{array}{cc}
I_2 & 0 \\
- \omega T^{-1} & 1
\end{array}
\right) \;. 
\label{T1}
\end{equation} 

Thus we arrive to the following \\ 
$\underline{Definition.}$ The dual to the Hopf superalgebra 
$osp_\xi(1|2)$ generated by 
the entries of $ {\bf T}$ (\ref{T}) subject to the relations 
(\ref{SL}) - (\ref{CO}) 
is said to be the quantum supergroup  $OSp_\xi(1|2)$. 

Another way to define this $OSp_\xi(1|2)$ is to use the twist 
element $\cal F$ as the pseudodifferential operator on 
the Lie supergroup $OSp(1|2)$, and redefine super-commutative 
product of functions on this supergroup. 

The reduction or Hopf superalgebra homomorphism, of $OSp_\xi(1|2)$ to  
$(sB_-)_\xi'$ is given by : 
$$b = \alpha = \beta = 0\,, \; g = 1 \,, \; 
a = d^{-1} = \exp (\nu )\,, \; 
\gamma \,a^{-1} = \delta = \frac {1}{2} \eta \,, \; c = 2 \xi x a \,. 
$$ 

\section{Conclusion} 
Using embedding of the Lie algebra $sl(2)$ into the rank one 
orthosymplectic superalgebra the latter one was deformed by the 
twist element ${\cal F} \in {\cal U}\bigl(sl(2)\bigr)^{\otimes 2}$. 
Although the deformed Lie superalgebra is finite dimensional 
it can be used for further deformation of infinite dimensional  
Hopf superalgebras (e.g. super-Yangians) and integrable 
models \cite{KS}. There are also possibilities for different 
contractions. The work in this direction is in progress. 

${\bf Acknowledgements}$ 
The authors are greatful to R. Kashaev and M. Rasetti for useful 
discussions. We appreciate the hospitality of the Institute for 
Scientific Interchange Foundation. This research was supported 
by the INTAS contract 94-1454 and by the RFFI grant 96-01-00851.

\end{document}